\documentclass[twoside,11pt,reqno]{amsart}
\usepackage{amsmath,amsthm,amssymb,amstext,amsfonts,amscd}
\usepackage{graphicx}
\usepackage{multirow}
\usepackage{url}
\usepackage{enumerate}
\begin{document}
\title[{\bf  A new proof of the fundamental two-term transformation  ...}]
{\bf A new proof of the fundamental two-term transformation for the series ${}_3F_2(1)$ due to Thomae }
\author[{\bf  A. K. Rathie}]{\bf Arjun K. Rathie}
\address{ Arjun K. Rathie, Department of Mathematics, Central University of Kerala, Kasaragad 671123, Kerala, India }
 \email{akrathie@cukerala.ac.in}

\begin{abstract}
The aim of this short note is to provide a very simple proof  for obtaining the fundamental two-term transformation for the series ${}_3F_2(1)$ due to Thomae.\\


\end{abstract}
  \maketitle
    \section{Introduction and Results Required} 
    
    We start with the following very useful and interesting Thomae[7] and Kummer[1] two-term transformations respectively as 
\begin{align}
 {}_3F_2  \left[\begin{array}{c} a,\; b, \; c \\ d, \; e \end{array} ; 1\right] & = \frac{\Gamma(d) \Gamma(e) \Gamma(d+e-a-b-c)}{\Gamma(a)\Gamma(d+e-a-b) \Gamma(d+e-a-c)}\;  \nonumber \\
         & \qquad . \;{}_3F_2  \left[\begin{array}{c} d-a,\; e-a, \;d+e-a-b-c \\ d+e-a-b,\; d+e-a-c \end{array} ; 1\right] 
 \end{align}     
 provided $\Re(a) > 0$ and $\Re(d+e-a-b-c) > 0$
and 
\begin{align}
 {}_3F_2  \left[\begin{array}{c} a,\; b, \; c \\ d,\; e \end{array} ; 1\right] 
& = \frac{ \Gamma(e) \Gamma(d+e-a-b-c)}{\Gamma(e-c)\Gamma(d+e-a-b) }\; \nonumber \\ 
   & \qquad .\; {}_3F_2  \left[\begin{array}{c} d-a,\; d-b, \;c \\d, \;  d+e-a-b\;  \end{array} ; 1\right] 
 \end{align}     
   provided $\Re(e-c) > 0$ and $\Re(d+e-a-b-c) > 0$.\\
   
As given in Bailey's tract[2], the Thomae transformation (1) can be established with the help of following classical Gauss's summation theorem[2, p.2, Eq. (1) or 5]
\begin{equation}
{}_2F_1 \left[\begin{array}{c} a,\; b \\ c \end{array}; 1\right] = \frac{\Gamma(c)\; \Gamma(c-a-b)}{\Gamma(c-a) \; \Gamma(c-b)}
\end{equation}
provided $\Re(c-a-b) > 0$.\\
 For a derivation of Kummer's transformation (2) we refer the standard text of Andrews, et al[1, eq. 3.3.5, p.142].
\\

In 1999, by considering the double series and summing up in two ways by employing Gauss's summation theorem (3) and following Saalsh\"utz summation theorem [5]
\begin{equation}
{}_3F_2\left[ \begin{array}{c} -n, \; a,\; b \\ c, \; 1+a+b-c-n \end{array}; 1\right] = \frac{(c-a)_n \; (c-b)_n}{(c)_n \; (c-a-b)_n}  
\end{equation} 

Exton[4] re-derived the Kummer's transformation (2).\\

In 2004, Rathie, et al.[6] have given a very short proof of Kummer's transformation (2). \\

The aim of this short note is to provide a very simple proof for obtaining the fundamental two-term transformation (1) for the series ${}_3F_2(1)$ due to Thomae.\\
For this, we require the following result known as Euler's second transformation[5] :
\begin{equation}
{}_2F_1 \left[\begin{array}{c}a,\; b \\ d \end{array}; x\right] = (1-x)^{d-a-b} 
\; {}_2F_1 \left[\begin{array}{c}d-a,\; d-b \\ d \end{array}; x\right] 
\end{equation}
\section{Derivation of Thomae Transformation (1) }
In order to derive Thomae transformation (1), we start with  the following known integral involving hypergeometric function[3, Equ.(5), p. 399 ] :
\begin{equation}
\int_0^1 x^{c-1} (1-x)^{e-c-1}\; {}_2F_1\left[\begin{array}{c}a,\; b \\ d \end{array}; x\right]dx = 
\frac{\Gamma(c) \Gamma(e-c)}{\Gamma(e)}\; {}_3F_2\left[ \begin{array}{c} a, \; b,\; c \\ d, \; e \end{array}; 1\right]
\end{equation}
provided $\Re(c)>0$, $\Re(e-c)>0$ and $\Re(d+e-a-b-c) > 0$\\

In order to prove (6), express ${}_2F_1$ as a series, change the order of integration and summation, which is justified  due to uniform convergence of the series, evaluate the integral with the help of the following beta integral :
\begin{equation}
\int_0^1 x^{\alpha-1} (1-x)^{\beta-1}  dx =  \frac{\Gamma(\alpha) \Gamma(\beta)}{\Gamma(\alpha+\beta)}
\end{equation}
provided $\Re(\alpha) > 0$ and $\Re(\beta) > 0$.\\
and then summing up the series, we easily arrive at the right-hand side of (6).\\

Now, in order to establish Thamae transformation (1), consider the result (6) in the form : 
\begin{align}
{}_3F_2 & \left[\begin{array}{c} a,\; b, \; c \\ d\; e \end{array} ; 1\right] \nonumber \\
& = \frac{ \Gamma(e) }{\Gamma(c) \; \Gamma(e-c)}\;\int_0^1 x^{c-1} (1-x)^{e-c-1} \; {}_2F_1\left[\begin{array}{c} a, b \\ d \end{array} ; x\right] dx  
\end{align}
Applying the result (5) to the ${}_2F_1$ on the right-hand side, we have
\begin{align}
{}_3F_2 & \left[\begin{array}{c} a,\; b, \; c \\ d, \; e \end{array} ; 1\right] \nonumber \\
& = \frac{ \Gamma(e) }{\Gamma(c) \; \Gamma(e-c)}\;\int_0^1 x^{c-1} (1-x)^{d+e-a-b-c-1} \; {}_2F_1\left[\begin{array}{c} d-a, d-b \\ d \end{array} ; x\right] dx  
\end{align}
Now, expressing ${}_2F_1$ as a series, change the order of integration and summation, which is easily seen to be justified due to uniform convergence of the series, we have 
\begin{align}
{}_3F_2 & \left[\begin{array}{c} a,\; b, \; c \\ d,\; e \end{array} ; 1\right] \nonumber \\
& = \frac{ \Gamma(e) }{\Gamma(c) \; \Gamma(e-c)}\; \sum_{n=0}^\infty \frac{(d-a)_n (d-b)_n}{(d)_n \; n!} . \int_0^1 x^{c+n-1} (1-x)^{d+e-a-b-c-1}  dx  
\end{align}
  Evaluating the integral with the help of (7), we have 
  \begin{align}
{}_3F_2 & \left[\begin{array}{c} a,\; b, \; c \\ d,\; e \end{array} ; 1\right] \nonumber \\
& = \frac{ \Gamma(e) \; \Gamma(d+e-a-b-c) }{\Gamma(e-c) \; \Gamma(d+e-a-b)}\; 
\sum_{n=0}^\infty \frac{(d-a)_n \; (d-b)_n\; (c)_n}{(d)_n \; (d+e-a-b)_n\; n!} . 
\end{align}
Summing up the series, we have
 \begin{align}
{}_3F_2 & \left[\begin{array}{c} a,\; b, \; c \\ d, \; e \end{array} ; 1\right] \nonumber \\
& = \frac{ \Gamma(e) \; \Gamma(d+e-a-b-c) }{\Gamma(e-c) \; \Gamma(d+e-a-b)}\; 
{}_3F_2  \left[\begin{array}{c} d-a,\; d-b, \; c \\ d,\; d+e-a-b \end{array} ; 1\right]
\end{align}
Now, writing (12) in the form
\begin{align}
{}_3F_2 & \left[\begin{array}{c} a,\; b, \; c \\ d,\; e \end{array} ; 1\right] \nonumber \\
& = \frac{ \Gamma(e) \; \Gamma(d+e-a-b-c) }{\Gamma(e-c) \; \Gamma(d+e-a-b)}\; 
{}_3F_2  \left[\begin{array}{c} d-b, \; c, \; d-a\\  d+e-a-b, \; d \end{array} ; 1\right]
\end{align}
Now using (6) to the ${}_3F_2$ on the right-hand side of (13), we have 
\begin{align}
{}_3F_2  \left[\begin{array}{c} a,\; b, \; c \\ d,\; e \end{array} ; 1\right] & = \frac{\Gamma(d)\; \Gamma(e) \; \Gamma(d+e-a-b-c) }{\Gamma(a)\; \Gamma(d-a)\; \Gamma(e-c) \; \Gamma(d+e-a-b)}\; \nonumber \\
& \qquad . \int_0^1  x^{d-a-1} (1-x)^{a-1} \; {}_2F_1\left[\begin{array}{c} d-b, \; c \\ d+e-a-b\end{array}; x \right] dx
\end{align}
Using Euler's second transformation, we have   
\begin{align}
{}_3F_2  \left[\begin{array}{c} a,\; b, \; c \\ d, \; e \end{array} ; 1\right] & = \frac{ \Gamma(d) \; \Gamma(e) \; \Gamma(d+e-a-b-c) }{ \Gamma(a)\; \Gamma(d-a)\;\Gamma(e-c) \; \Gamma(d+e-a-b)}\; \nonumber \\
& \qquad .\int_0^1  x^{d-a-1} (1-x)^{e-c-1} \; {}_2F_1\left[\begin{array}{c} e-a, d+e-a-b-c \;  \\ d+e-a-b\end{array}; x \right]dx
\end{align}
Finally, expressing ${}_2F_1$ as a series, change the order of integration and summation and after some algebra, evaluating the beta integral, we have
\begin{align}
{}_3F_2  \left[\begin{array}{c} a,\; b, \; c \\ d,\; e \end{array} ; 1\right] & = \frac{\Gamma(d) \; \Gamma(e) \; \Gamma(d+e-a-b-c) }{ \Gamma(a)\; \Gamma(d+e-a-b)\;\Gamma(d+e-a-c) }\; \nonumber \\
&  \qquad . \sum_{n=0}^\infty \frac{(d-a)_n \; (e-a)_n \; (d+e-a-b-c)_n}{(d+e-a-b)_n \; (d+e-a-c)_n \; n!}
\end{align}
Finally, summing up the series, we get the required Thomae transformation formula (1)\\
This completes the proof of the   Thomae transformation formula.

\end{document}